\documentclass[12pt,a4paper,twoside]{article}
\usepackage{amsmath}
\usepackage{amsfonts}
\usepackage{amssymb}
\usepackage[left=2.00cm, right=2.00cm, top=3.00cm, bottom=3.00cm]{geometry}
\usepackage{hyperref}
\usepackage{cases}
\usepackage{graphicx}
\usepackage{xypic}
\usepackage{theorem}
\hypersetup{ colorlinks=true
pdfborder=0 0 0
}
\numberwithin{equation}{section}
\newcommand{\cc}[1]{\mathcal{C}_{#1}^{+}}
\newcommand{\ccc}[1]{\mathcal{C}_{#1}^{-}}

\newcommand{\p}[1]{P_{#1}}
\newcommand{\pp}[1]{P^*_{#1}}
\newcommand{\ty}[1]{\widetilde{T}_{#1}}
\newcommand{\con}[2]{\substack{#1\leq #2w_v\\l(#1w_v)=l(#1)+l(w_v)}}

\def\d{\mathcal{D}_0}

\def\qed{\hfill$\square$}

\def\ra{\rightarrow}
\def\irr{\operatorname{Irr}}\def\h{\mathcal{H}}\def\zg{\mathbb{Z}[\Gamma]}
\def\z{\mathcal{Z}_\mathbb{Z}}
\def\mspec{\operatorname{MaxSpec}(\mathcal{Z}_k)}\def\m{\mathfrak{m}}
\def\ezd{\operatorname{End}_{\mathcal{Z}}(\Delta)}
\def\jg{J_{0,\Gamma}}

\def\pr {\mathfrak{p}}
\newtheorem{thm}{Theorem}[section]
\newtheorem{prop}[thm]{Proposition}
\newtheorem{lemma}[thm]{Lemma}
\newtheorem{cor}[thm]{Corollary}
\newtheorem{defn}[thm]{Definition}
\newtheorem{assumption}[thm]{Assumption}

\newtheorem{remark}[thm]{Remark}
\newtheorem{thm_def}[thm]{Theorem (Definition)}
\theoremstyle{definition}
\author{}
\title{The based ring the lowest generalized two-sided cell\\
of an extended affine Weyl group}

\begin{document}
\author{Xun Xie} 
\date{}
\maketitle
\allowdisplaybreaks
\bibliographystyle{alpha}
\abstract{Let $\mathbf{c}_0$  be the lowest generalized two-sided cell of an extended affine Weyl group W. We determine the structure of the based ring of $\mathbf{c}_0$. For this we show that certain conjectures of Lusztig on generalized cells (called P1-P15) hold for $\mathbf{c}_0$. As an application, we use the structure of the based ring to study certain simple modules of Hecke algebras of $ W $ with unequal parameters, namely those attached to $\mathbf{c}_0$.
	Also we give a  set of prime ideals $\mathfrak{p}$ of the center $\mathcal{Z}$  of the generic affine Hecke algebra $\mathcal{H}$ such that the reduced affine Hecke algebra $k_\mathfrak{p}\mathcal{H}$ is  simple over $k_\mathfrak{p}$, where $k_\mathfrak{p}=\operatorname{Frac}(\mathcal{Z}/\pr)$ is the residue field of $\mathcal{Z}$ at $\mathfrak{p}$. In particular, we show that the algebra $\mathcal{H}\otimes_\mathcal{Z}\operatorname{Frac}(\mathcal{Z})$ is a split  simple algebra over the field $ \operatorname{Frac}(\mathcal{Z})$.}

\section{Introduction}

We are concerned with the based ring (also called asymptotic ring) of the lowest generalized two-sided cell of an extended affine Weyl group $ W $. The based ring of an two-sided cell of a certain Coxeter group was introduced by Lusztig in \cite{lusztig1987cellsII}, which is useful in studying representations of Hecke algebras (see for examples \cite{lusztig1987cellsIII,lusztig1989cells,Xi1990based,xi1994representations,xi2007representations}).

Generalized cells of Coxeter groups are defined by Lusztig using Hecke algebras of Coxeter groups with unequal parameters in  \cite{lusztig1983left}. For a certain Coxeter group, the based ring of a generalized two-sided cell is defined in \cite[\S18]{lusztig2003hecke}, by assuming 15 conjectural properties of the generic Hecke algebras (or generalized cells) of the Coxeter group with unequal parameters,  which are called properties P1-P15 (see \cite[\S14]{lusztig2003hecke}).

In this paper, we show that the properties P1-P15 are valid for the lowest generalized two-sided cell of an extended affine Weyl group. This implies that the based ring of the lowest generalized two-sided cell is well defined. We then determine the structure of the based ring. As an application, we use these results to study representations of affine Hecke algebras with unequal parameters attached to the lowest generalized two-sided cell.

The contents of the paper are as follows. In section 2 we recall some basic facts about extended affine Weyl groups and their Hecke algebras and the lowest generalized two-sided cell of an extended affine Weyl group. In section 3 we study a formula of Xi, which gives a decomposition of the (generalized) Kazhdan-Lusztig element $C_w$ for $w$ being in the lowest generalized two-sided cell. This formula is a key for determining the structure of the based ring of the lowest generalized two-sided cell. In section 4 we show that the properties P1-P15 are valid for  the lowest two-sided cell of an extended affine Weyl group. In section 5 we determine the based ring of the lowest generalized two-sided cell. In section 6 we consider the based ring of the lowest generalized two-sided cell of an affine Weyl group of type $\tilde C_n$. In section 7 we use the ideas in \cite{xi1994representations,xi2007representations} to study simple representations of an affine Hecke algebras with unequal parameters attached to the lowest generalized two-sided cell. In section 8 we study the reduced affine Hecke algebra $k_\mathfrak{p}\mathcal{H}=\mathcal{H}\otimes_\mathcal{Z} k_\mathfrak{p}$, where $\mathcal{H}$ is the generic affine Hecke algebra, $\mathcal{Z}$ is the center of $\mathcal{H}$, $\mathfrak{p}$  is a prime ideal of $\mathcal{Z}$, $k_\mathfrak{p}=\operatorname{Frac}(\mathcal{Z}/\pr)$ is the residue field of $\mathcal{Z}$ at $\mathfrak{p}$.

\section{Preliminaries}\label{pre}

In this section, we recall some basic definitions and  facts about affine Weyl groups and their Hecke algebras and the lowest generalized two-sided cell of an extended affine Weyl group.

\subsection{Extended affine Weyl groups}\label{extended aff}

Let $R$ be an irreducible  root system. Its Weyl group, root lattice and weight lattice are denoted by $W_0$, $Q$ and $P$ repectively. Then  $W'=W_0\ltimes Q$ is an affine Weyl group and $W=W_0\ltimes P= \Omega\ltimes W'$  is an  extended affine Weyl group, where $\Omega$ is a finite subgroup of $W$ and is isomorphic to $P/Q$. The length function $l:W'\to \mathbb{N}$ can be extended to a length function $l:W\to \mathbb{N}$ by defining $l(\omega w)=l(w)$ for $\omega\in\Omega$ and $w\in W'$. Denote by $e$ the neutral element of the group $W$.

Let $E=\mathbb{R}\otimes Q$ and $(~,~)$  be an $W_0$-invariant inner product on $E$. Let $\mathcal{F}$ be the set of hyperplanes in $E$ consisting of all hyperplanes
$$H_{\alpha,n}=\{x\in E\,|\, \frac{(x,2\alpha)}{(\alpha,\alpha)}=n\},\qquad \alpha\in R,\ \ n\in\mathbb{Z}.$$

Let $X$ be the set of alcoves, i.e. the set of connected components of the set $E-\bigcup_{H\in \mathcal{F}}H$. Then $W'$ is isomorphic to the group generated by all reflections on $E$ with respect to the hyperplanes in $\mathcal{F}$. Thus $W'$ acts on $E$ and $X$ naturally.  We regard this action as right action. Another action of $W'$ on $X$, regarded as a left action, was introduced in \cite[1.1]{Lusztig1980Janzten}. The left action and the right action commute. For every simple reflection $s\in W'$ and every alcove $A\in X$, $A$ and $sA$ has a common face and the common face of $A$ and $sA$ is called a face (of $A$) of type $s$.

Let $\Gamma$ be an abelian group written additively.
A \textit{weight function} of the {extended} affine Weyl group, $L:W\longrightarrow\Gamma$, is defined by the condition: $L(ww')=L(w)+L(w')$ if $l(ww')=l(w)+l(w')$, and $L(\pi)=0\in \Gamma$ if $\pi\in \Omega$. Denote by $\mathbb{Z}[\Gamma]$ the group ring of the group $\Gamma$. Let $q$ be a symbol. Following \cite{bonnafe2006tw}, an element of $\mathbb{Z}[\Gamma]$ is written in the form $\sum_{\gamma\in\Gamma}a_\gamma q^\gamma$ with $a_\gamma\in \mathbb{Z}$. We often write $q_w$ for $q^{L(w)}$. In particular, $q_\pi=q^0=1\in q^{\Gamma}$ for each $\pi\in \Omega$.

Given a weight function $L$ on $W$, we can associate with each hyperplane $H\in \mathcal{F}$ an element $L(H)=L(s)$ in $\Gamma$, where $s$ is an simple reflection such that $H$ supports a face of type $s$. This is well defined by \cite[Lemma 2.1]{bremke1997generalized}. Furthermore, by \cite[Lemma 2.2]{bremke1997generalized}, we have
\begin{itemize}
\item[(a)] $L(H)=L(H')$ when $H$ and $H'$ are parallel.
\end{itemize}

\begin{remark}A weight function on $W'$ may not be extended to a weight function on $W$. The property (a) does not  hold for some  weight functions of $W'$ of type $\tilde C_r, r\geq1$,  see Section \ref{nonext}.
\end{remark}

\subsection{Affine Hecke algebras}\label{aff hec alg}

\begin{defn}
Let $W$ be an extended affine Weyl group and let $L:W\longrightarrow\Gamma$ be a weight function of $W$. The affine Hecke algebra $\mathcal{H}$ with respect to weight function $L$ is defined to be  the  algebra over $\mathbb{Z}[\Gamma]$ with generators $\{T_w\mid w\in W\}$ and relations
\begin{itemize}
\item[(i)] $T_wT_{w'}=T_{ww'}$ if $l(ww')=l(w)+l(w')$,
\item[(ii)] $(T_s+1)(T_s-q_s ^2)=0$ if $s\in S$.
\end{itemize}
We often use the notation $\ty{w}:=q_w^{-1}T_w$, $w\in W$. The subalgebra generated by $T_w,w\in W'$ is denoted by $\mathcal{H'}$.
\end{defn}

Assume that $\Gamma$ is a totally ordered abelian group. Then there is a total order $\leq$ on $\Gamma$ such that $\gamma_1+\gamma<\gamma_2+\gamma$ for any $\gamma\in \Gamma$ whenever $\gamma_1<\gamma_2$. Now the (generalized) Kazhdan-Lusztig basis (KL-basis) $\{C_w \mid w\in W  \}$ of $\mathcal{H}$ with respect to the pair $(L,\leq)$ is defined by the following property:
\begin{itemize}
\item[(a)]
\begin{itemize}
\item $\overline{C}_w=C_w$, where $\overline{\,\cdot\,}$ is the involution of algebras on $\mathcal{H}$ such that $\overline{q^\gamma}=q^{-\gamma}$, $\overline{T_s}=T_s^{-1}$ and $\overline{T_\pi}=T_\pi$, for $\gamma\in\Gamma$, $s\in S$, $\pi\in \Omega$.
\item $C_w\equiv \ty{w}\mod{\mathbb{Z}[\Gamma^{<0}]\{\ty{w}|w\in W\}}$, where $\mathbb{Z}[\Gamma^{<0}]$ is the set of $\mathbb{Z}$-linear combinations of the elements $\gamma$ in $\Gamma$ such that  $\gamma<0$.
\end{itemize}
\end{itemize}

The element $C_w$ can be expressed  as $C_w=\sum_{y\leq w}\pp{y,w}\ty{y}=q_w^{-1}\sum_{y\leq w}\p{y,w}T_y$ where $\p{y,w}$ is a polynomial in $q_s^2,s\in S$, and $\pp{y,w} $ is a polynomial in $q_s^{-1},s\in S$.

As in \cite{lusztig2003hecke}, using the KL-basis, one can  define preorders $\leq_L$, $\leq_R$, $\leq_{LR}$ on $W$.  These preorders induce naturally equivalence relations $\sim_L$, $\sim_R$, $\sim_{LR}$ on $W$. The corresponding equivalence classes are called respectively generalized left cells, generalized right cells, and generalized two-sided cells of $W$ (with respect to the pair $(L,\leq)$).

\begin{itemize}
\item[(b)] Since $\Gamma$ is totally ordered, we can define the degree map $$\deg:\mathbb{Z}[\Gamma]\longrightarrow\Gamma$$ by setting $\deg(\sum_{\gamma\in \Gamma}a_\gamma q^{\gamma})=\max\{\gamma|a_\gamma\neq0\}$.
\end{itemize}

The elements $h_{x,y,z}$, $m_{x,y,z}$ of $\mathbb{Z}[\Gamma]$ are defined respectively by the formulas $$C_xC_y=\sum_{z}h_{x,y,z}C_z,\quad \ty{x}\ty{y}=\sum_{z}m_{x,y,z}\ty{z}.$$Then\begin{itemize}\item [(c)]
 Lusztig's \textbf{a}-function is defined to be $$\textbf{a}(z)=\sup\{\deg(h_{x,y,z})~|~x,y\in W \}\in\Gamma\cup\{\infty\}.$$
 \end{itemize}

  If $\textbf{a}(z)\neq\infty$, then $h_{x,y,z}$ can be written in the form $$h_{x,y,z}=\gamma_{x,y,z^{-1}}q^{\textbf{a}(z)} \text{+ lower degree terms},$$ where $\gamma_{x,y,z^{-1}}\in \mathbb{Z}$. In this case there exists some $x,y\in W$ such that $\gamma_{x,y,z^{-1}}$ is nonzero.

\begin{assumption}\label{assum1}
Hereafter, we will always assume that $\Gamma$ is a totally ordered abelian group and $L$ is a positive weight function, i.e. $L(s)>0$ for all $s\in S$.
\end{assumption}

\subsection{Lowest two-sided cells}\label{low_c_0}
In this subsection we recall some known results on the lowest generalized two-sided cell of the extended affine Weyl group $W$.

The lowest two-sided cell of $W$ is described in \cite{shi1987two,shi1988two}. The lowest generalized two-sided cell of $W$ is described in \cite[Theorem 3.22]{xi1994representations}. It turns out the generalized lowest two-sided and lowest two-sided cell of $W$ coincide. It is
$$\mathbf{c}_0=\{ww_0u\in W\ |\ l(ww_0u)=l(w)+l(w_0)+l(u)\},$$
where $w_0$ is the longest element of $W_0$. We shall keep to use the term generalized lowest two-sided cell for emphasizing its connection with affine Hecke algebras with unequal parameters.

Following \cite{shi1988two} we give a more detailed description of $\mathbf{c}_0$.

Keep the notations in \ref{extended aff}. Let $v$ be a special point in $E$ (that is, there are $|R|/2$ hyperplanes in $\mathcal{F}$ passing through $v$) and let $W_v$ be the stabilizer in $W'$ (under left action) of the set of alcoves containing $v$ in its closure. Let $w_v$ be the longest element of $W_v$. A connected component of $E-\cup_{v\in H}H$ is called a quarter with vertex $v$. For every special point $v$, we will fix a quarter $\cc{v}$ such that for any two special points $v, v'$, $\cc{v'}$ is a translate of $\cc{v}$.  Denote by $A_v^+$ the unique alcove contained in $\cc{v}$ and containing $v$ in its closure. Let $\ccc{v}=-\cc{v}$ and $A_v^-=w_vA_v^+\subset\ccc{v}$. Let $\mathcal{F}^*$ be the set of hyperplanes $H\in \mathcal{F}$ such that $H$ is parallel to a wall of $\cc{v}$. A connected component of $E-\cup_{H\in \mathcal{F}^*}H$ is called a box. The box containing $A_v^+$ is unique and is denoted by $\Pi_v^+$. For $H\in\mathcal{F}$, let $E_H^+$ be the connected component of $E-H$ that meets $\cc{v}$ for any special point $v$.

Let $$B_v=\{\pi w\in W\vert w\in W', wA_v^+\subset \Pi_v^+,\pi\in \Omega\},$$ $$U_v=\{\pi w\in W\vert w\in W', wA_v^+\subset \cc{v},\pi\in \Omega\}.$$ The lowest two-sided cell of $W$ with respect to $(L,\leq)$ can be described as (see [Shi88] and \cite[Thm. 6.13]{bremke1997generalized}).

\begin{itemize}
\item[(a)]
$\mathbf{c}_0=\{ww_vw'^{-1}~|~w\in B_v,w'\in U_v\}=\{z\in W~|~\textbf{a}(z)=L(w_0) \}$.
\end{itemize}

The cell $\mathbf{c}_0$ is independent of the choice of the special point $v$, since if  $v,v'$ are two special points, then there always exists $\pi\in \Omega$ such that $\pi w_v\pi^{-1}=w_{v'}$. Recall that $\Omega$ is the finite subgroup of $W$ consisting of elements of length 0 and that $\Omega$ normalizes $W'$.  The lowest two-sided cell $\mathbf{c}_0$ can  also be described as $\mathbf{c}_0=\{ww_vw'^{-1}~|~w\in U_v,w'\in B_v\}$. Note that we always have $l(ww_vw'^{-1})=l(w)+l(w_0)+l(w')$ for $w\in B_v$, $w'\in U_{v}$ or $w\in U_v,w'\in B_v$. It is known that the set
$$\mathcal{H}^{\mathbf{c}_0}=\mathbb{Z}[\Gamma]\{C_w \mid w\in\mathbf{c}_0\}$$ is a two-sided ideal of the affine Hecke algebra $\mathcal{H}$.

For $x\in P$, we denote by $p_x$ the corresponding element in $W=W_0\ltimes P$. The set of dominant weights is defined to be $P^+=\{x\in P~|~l(w_0p_x)=l(w_0)+l(p_x) \} $. We have

\begin{itemize}
\item[(b)] Each $z\in \mathbf{c}_0$ can be written uniquely as $w_1w_0p_xw_2^{-1}$ with $w_1,w_2\in B_0$, $x\in P^+$, see \cite[3.2]{Xi1990based}.
\end{itemize}

 For every $x\in P$, we can find $y,z\in P^+$ such that $x=y-z$ and set  $\theta_x=\ty{p_y}(\ty{p_z})^{-1}$.  It is known that $ \theta_x $ is well defined and is independent of the choice of $ y,z $. For $x\in P^+$, we define $S_x=\sum_{x'\in P}d(x',x)\theta_{x'}$, where $d(x',x)$ is the dimension of the weight space $V(x)_{x'}$, and  $V(x)$ is an irreducible module of highest weight $x$ of the  simply connected complex algebraic group $G_\mathbb{C}$ with Weyl group $W_0$. Then
\begin{itemize}
\item[(d)] The elements $S_x$, $x\in P^+$ form a $\mathbb{Z}[\Gamma]$-basis of the center $\mathcal{Z}$ of $\mathcal H$. And $S_{x}S_{x'}=\sum_{x''\in P^+}m(x,x',x'')S_{x''}$, where $m(x,x',x'')$ is the multiplicity of the irreduicble module $V(x'')$ in the tensor product $V(x)\otimes V(x')$ (see \cite[2.9]{xi1994representations}).
\item[(e)] $C_{w_0p_x}=C_{w_0}S_x$ for $x\in P^+$. Hence $C_{w_0p_x}C_{w_0p_y}=\sum_{z\in P^+}m(x,y,z)h_{w_0,w_0,w_0}C_{w_0z}$ for $x,y\in P^+$(see \cite[3.21]{xi1994representations}).
\end{itemize}

\begin{defn}\label{center form}
Define $\mathcal{Z}_{\mathbb{Z}}$ to be the $\mathbb{Z}$-submodule of $\mathcal{Z}$ generated by $\{S_x\,|\,x\in P^+ \}$.
\end{defn}

By (d), $\mathcal{Z}_{\mathbb{Z}}$ is a $\mathbb{Z}$-lattice of $\mathcal{Z}$ (that is $\mathcal{Z}=\mathbb{Z}[\Gamma]\otimes_\mathbb{Z}\mathcal{Z}_{\mathbb{Z}}$), and $\mathcal{Z}_{\mathbb{Z}}$ is  isomorphic to the representation ring of the algebraic group $G_\mathbb{C}$.

The following facts are  well known or easy to check and will be used latter.
\begin{lemma}\label{predecomp}Let $\mathcal{H}$ be the Hecke algebra of $W$ with respect to a positive weight function $L:W\longrightarrow \Gamma$. Then
\begin{itemize}
\item[(i)]$\pp{u,w}=q_s^{-1}\pp{us,w}$ if $u<us\leq w$, $ws<w$. In particular, $\pp{y,w_v}=q_{yw_v}^{-1}=q^{L(y)-L(w_v)}$, for $y\leq w_v$. Hence $C_{w_v}=\sum_{y\leq w_v}q_{yw_v}^{-1}\ty{y}=q_{w_v}^{-1}\sum_{y\leq w_v} T_y$.

\item[(ii)] $h_{w_v,w_v,w_v}$ is equal to ${q_{w_v}^{-1}}\sum_{y\in W_v} q_y^2$ and is nonzero in $\mathbb{Z}[\Gamma]$ since $\Gamma$ is totally ordered.
\item[(iii)] $C_{ww_v}=E_wC_{w_v}$ for $w\in U_v$ where $E_w=\sum_{\substack{u\leq ww_v\\l(uw_v)=l(u)+l(w_v)}}\pp{uw_v,ww_v}\ty{u}$.

\item[(iv)] $C_{w_vw^{-1}}=C_{w_v}F_w$ for $w\in U_v$ where $F_w=\sum_{\con{u}{w}}\pp{uw_v,ww_v}\ty{u^{-1}}$.

\item[(v)] $\{{ww_v}\vert w\in U_v\}$ is a left cell of $W$. 
\end{itemize}
\end{lemma}
\textit{Proof.} (i) and (ii) are well known, see for example \cite[1.14.13 and 1.17(i)(ii)]{xi1994representations}. (iii) can be calculated straightforward as follows.
\begin{align*}
C_{ww_v}&=\sum_{z\leq ww_v}\pp{z,ww_v}\ty{z}\\
&=\sum_{\substack{u\leq ww_v\\l(uw_v)=l(u)+l(w_v)\\y\leq w_v}}\pp{uy,ww_v}\ty{u}\ty{y}\\
&=\sum_{\substack{u\leq ww_v\\l(uw_v)=l(u)+l(w_v)\\y\leq w_v}}q_{yw_v}^{-1}\pp{uw_v,ww_v}\ty{u}\ty{y}& \text{(by (i))}\\
&=\left(\sum_{\con{u}{w}}\pp{uw_v,ww_v}\ty{u}\right)\left(\sum_{y\leq w_v}q_{yw_v}^{-1}\ty{y}\right)\\&=E_wC_{w_v}.& \text{(by (i))}
\end{align*}

(iv) Let ${\cdot\,}^\flat:\mathcal{H}\longrightarrow \mathcal{H}$ be the anti-isomorphism of algebras over $\mathbb{Z}[\Gamma]$ such that $(\ty{u})^\flat=\ty{u^{-1}}$ for $u\in W$. Then ${\cdot\,}^{\flat}$ commutes with the bar involution $\overline{\,\cdot\,}$ in \ref{aff hec alg}(a), and hence $(C_w)^\flat=C_{w^{-1}}$. Therefore $C_{w_vw^{-1}}=(C_{ww_v})^\flat=(E_wC_{w_v})^\flat=C_{w_v}F_w$, where we use the fact that $w_v^2=1$ and $(E_w)^\flat=F_w$.

(v) follows from the proof of Theorem 3.22 in \cite{xi1994representations}. \qed

\section{A formula of Xi}\label{formula}

In \cite[Theorem 2.9]{xi1994based}, Xi established a decomposition formula for $C_w, \ w\in \mathbf{c}_0$ when the weight function is constant. The formula is reproved in \cite{Blasiak2009fac}. In this section we will show that Xi's formula remains valid for general (positive) weight functions.

For $A$, $B\in X$ are two alcoves, the distant function $d(A,B)$ is  defined to be the number of hyperplanes separating $A$ from $B$ counting by signs (see \cite[1.4]{Lusztig1980Janzten}), and there exists a partial order $<$ on $X$ induced by the alcoves with distant 1, \cite[1.5]{Lusztig1980Janzten}. The following technical lemma is useful, see \cite[Lemma 4.3]{Lusztig1980Janzten} and \cite[Lemma 4.3]{bremke1997generalized}.

\begin{lemma}[Lusztig, Bremke]\label{lus_lemma}Let $W$ be an extended affine Weyl group with a positive weight fuction $L$.
Let $v$ be a special point, $A$ be an alcove containing $v$ in its closure and $y\in W_v$ such that $A=yA_v^+$. Let $s_1$,$\dots$,$s_k\in S$ be such that $d(A_v^+,s_k\cdots s_1A_v^+)=k$ and let $1\leq i_1,\dots i_p\leq k$ be such that \begin{gather*}s_{i_t}\cdots\hat{s}_{i_{t-1}}\cdots\hat{s}_{i_1}\cdots s_1(A)<\hat{s}_{i_t}\cdots\hat{s}_{i_{t-1}}\cdots\hat{s}_{i_1}\cdots s_1(A)\end{gather*}
for $t=1,\cdots, p$. Here $\hat{s}_{i_t}$ means deleting $s_{i_t}$. Then we have
\begin{itemize}
\item[(i)] $L({s_{i_1}})+\cdots+ L({s_{i_p}})\leq L(y)$.
\item[(ii)] If moreover $s_k\cdots s_1(A_v^+)\subset \Pi_v^+$ and $A\neq A_v^+$ (i.e. $y\neq1$), then the strict inequality holds, i.e., $L({s_{i_1}})+\cdots+ L({s_{i_p}})< L(y)$.
\end{itemize}\qed
\end{lemma}

Let $\mathcal{M}$ be the $\mathbb{Z}[\Gamma]$-module with basis $X$. Define $\mathcal{M}$ to be an $\mathcal{H'}$-module via
\begin{numcases}{\ty{s}{A}=}\label{tma}
{sA} & if $sA>A$\\ \notag
{sA}+(q_{s}-q_s^{-1}){A}& if $sA<A$
\end{numcases}
for $s\in S$.

For $w\in W'$, $A\in X$, if we write
\[\ty{w}{A}=\sum_{B\in X}\pi_{w,A,B}{B},\quad \pi_{w,A,B}\in\mathbb{Z}[\Gamma],\]
then $\pi_{w,A,B}$ actually are polynomials in $\xi_s=q_s-q_s^{-1}, s\in S$ with integral and positive coefficients.

The following corollary of the above lemma originates from \cite[4.2]{Lusztig1980Janzten}.
\begin{cor}\label{estimate}
Let $v$ be a special point, $u\in W'$ be such that $uA_v^+\subset\Pi_v^+$ and $y\in W_v $ be such that $A=yA_v^+\neq A_v^+$.  Recall that the degree map is defined by \ref{aff hec alg}(b). Then we have \[\deg(\pi_{u,A,B})<L({y})\]for any $B\in X$, where $\pi_{u,A,B}$ is the coefficients of ${B}$ in $\ty{u}{A}$, i.e. $\ty{u}{A}=\sum_{B\in X}\pi_{u,A,B}{B}$.
\end{cor}
\textit{Proof.} Assume that $u$ has a reduced expression $u=s_k\cdots s_1$ with $s_i\in S$, $1\leq i\leq k$. Then $uA_v^ +\subset\Pi_v^+ $ implies that $d(A_v^+,s_k\cdots s_1A_v^+)=k$.

By definition, we have \[\ty{u}{A}=\ty{s_k}\cdots \ty{s_1}(A)=\sum_{I\in\mathcal{I}}\prod_{k=1}^{p_I}(q_{s_{i_k}}-q_{s_{i_k}}^{-1}){s_k\cdots \hat{s}_{i_{p_I}}\cdots \hat{s}_{i_2}\cdots \hat{s}_{i_1}\cdots s_1(A)}\] where $\mathcal{I}$ is the set of all the sequences $(i_1<i_2<\cdots <i_{p_I})$ such that  \[s_t\cdots \hat{s}_{i_{t-1}}\cdots \hat{s}_{i_2}\cdots \hat{s}_{i_1}\cdots s_1(A)<\hat{s}_t\cdots \hat{s}_{i_{t-1}}\cdots \hat{s}_{i_2}\cdots \hat{s}_{i_1}\cdots s_1(A)\] for $t=1,\dots,p_I$. Now using Lemma \autoref{lus_lemma} (ii), we obtain the conclusion.\qed

\begin{cor}\label{est2}
Let $v$ be a special point, $u\in W'$ be such that $uA_v^+\subset\Pi_v^+$, $u'\in W'$ be such that $u' A_v^+\subset \cc{v}$ and $y\in W_v$ be such that $y<w_v$. Write the product $\ty{u}\ty{y}\ty{{u'}^{-1}}$ as $\sum_{z}a_{u,y,u'}^z\ty{z}$ with $a_{u,y,u'}^z\in\mathbb{Z}[\Gamma]$. Then we have $\deg(a_{u,y,u'}^z)<L({yw_v})$ in $\Gamma$.
\end{cor}

\textit{Proof.} The idea of the proof is same as that in \cite[7.9]{Lusztig1985cells} whose aim is to prove the boundedness of $\textbf{a}$-function for an affine Weyl group.

Let $A=u'A_v^-$. Then $A\subset \mathcal{C}_v^-$. Hence $\ty{{u'}^{-1}}{A}={A_v^-}$. This implies that $\ty{y}\ty{{u'}^{-1}}{A}={yA_v^{-}}={yw_vA_v^+}$.

On one hand, $$\ty{u}\ty{y}\ty{{u'}^{-1}}{A}=\ty{u}{yw_vA_v^+}=\sum_B\pi_{u,yw_vA_v^+,B}{B}.$$ On the other hand, $$\ty{u}\ty{y}\ty{{u'}^{-1}}{A}=\sum_z a_{u,y,u'}^z\ty{z}{A}=\sum_{z,B}a_{u,y,u'}^z\pi_{z,A,B}{B}.$$ Thus we have $$\sum_{z}a_{u,y,u'}^z\pi_{z,A,B}=\pi_{u,yw_vA_v^+,B}.$$ Since $a_{u,y,u'}^z$, $\pi_{z,A,B}$, and $\pi_{u,yw_vA_v^+,B}$ are all polynomials in $\xi_s$ with positive integral coefficients, we have $$\deg(a_{u,y,u'}^z\pi_{z,A,B})\leq \deg(\pi_{u,yw_vA_v^+,B}).$$ By formula \eqref{tma}, we know that $\pi_{z,A,zA} $ has constant term 1. Therefore, for $B=zA$, $$\deg(a_{u,y,u'}^z)\leq\deg(\sum_{z}a_{u,y,u'}^z\pi_{z,A,B})\leq\deg(\pi_{u,yw_vA_v^+,B})<L({yw_v}).$$ The last strict inequality is deduced from  Corollary \ref{estimate}.\qed

\begin{thm}\label{J0}
Let $v$ be a special point, $w\in W'$ be such that $wA_v^+\subset\Pi_v^+$, $w'\in W'$ be such that $w'A_v^+\subset\cc{v}$. Then\begin{gather*}
C_{ww_v{w'}^{-1}}=E_wC_{w_v}F_{w'}
\end{gather*}
where $E_w=\sum_{\substack{u\leq ww_v\\l(uw_v)=l(u)+l(w_v)}}\pp{uw_v,ww_v}\ty{u}$, $F_w=\sum_{\con{u}{w}}\pp{uw_v,ww_v}\ty{u^{-1}}$. 
\end{thm}
\textit{Proof.} On one hand,\begin{equation*}\begin{aligned}
E_wC_{w_v}F_{w'}&=\sum_{u,y,u'}\pp{uw_v,ww_v}q_{yw_v}^{-1}\pp{u'w_v,w'w_v}\ty{u}\ty{y}\ty{{u'}^{-1}}\\
&=\sum_{u,y,u',z}\pp{uw_v,ww_v}\pp{u'w_v,w'w_v}(q_{yw_v}^{-1}a_{u,y,u'}^z)\ty{z},
\end{aligned}\end{equation*}
where $u$ runs over $u\leq ww_v$, $l(uw_v)=l(u)+l(w_v)$, $u'$  runs over $u'\leq w'w_v$, $l(u'w_v)=l(u')+l(w_v)$, and $y$ runs over $y\leq w_v$. The element $\pp{uw_v,ww_v}\pp{u'w_v,w'w_v}$ has negative degree unless $u=w$ and $u'=w'$. And by Corollary \ref{est2}, $q_{yw_v}^{-1}a_{u,y,u'}^z$ has negative degree unless $y=w_v$. Thus $$E_wC_{w_v}F_{w'}\equiv \ty{w}\ty{w_v}\ty{w'^{-1}}\mod{\mathbb{Z}[\Gamma^{<0}]}.$$

On the other hand, $$E_wC_{w_v}F_{w'}=\frac1{h_{w_v,w_v,w_v}}C_{ww_v}C_{w_vw'^{-1}},\quad \overline{{h}_{w_v,w_v,w_v}}=h_{w_v,w_v,w_v}$$ (see Lemma \ref{predecomp}(iii)(iv)), hence $E_wC_{w_v}F_{w'}$ is $\overline{\,\cdot\,}$-invariant. The theorem  then follows from the definition of KL-basis (\ref{aff hec alg}(a)).\qed
\begin{thm}\label{uXi}
Let $v,v'$ be two special points. Let $u\in W'$ be such that $uA_v^+=A_{v'}^+\subset\cc{v}$, let $w\in W'$ be such that $wA_{v'}^+\subset \Pi_{v'}^+$, and let $w'\in W'$ be such that $w'A_{v}^+\subset \Pi_{v}^+$. Then
\[
C_{wuw_vw'^{-1}}=E_wC_{uw_v}F_{w'}.
\]
\end{thm}
\textit{Proof.} 
Since $u$ is a translation such that $uA_v^+=A_{v'}^+$, we can write $uw_v=w_{v'}u'^{-1}$ for some $u'$. Then $u'=w_vu^{-1}w_{v'}$ and $u'A_{v'}^-=w_vu^{-1}A_{v'}^+=w_vA_v^+=A_v^-\subset\mathcal{C}_{v'}^-$. This implies that $u'A_{v'}^+\subset \mathcal{C}_{v'}^+$. Using the theorem above for the special point $v'$, we get $$E_wC_{uw_v}=E_wC_{w_{v'}u'^{-1}}=E_wC_{w_{v'}}F_{u'}=C_{ww_{v'}u'^{-1}}=C_{wuw_v}.$$ Using the theorem above for the special point $v$ we see $C_{wuw_v}=E_{wu}C_{w_v}$. Again using the theorem above for the special point $v$, we get $C_{wuw_vw'^{-1}}=E_{wu}C_{w_v}F_{w'}=E_wC_{uw_v}F_{w'}$.\qed

Recall from \ref{low_c_0}(b) that, each $z\in \mathbf{c}_0$ has the form $w_1w_0p_xw_2^{-1}$ with $w_1,w_2\in B_0$, $x\in P^+$. The following  decomposition formula for the element $C_{w_1w_0p_xw_2^{-1}}$ is useful.
\begin{thm}[Xi's formula]\label{0xi}
$C_{w_1w_0p_xw_2^{-1}}=E_{w_1}C_{w_0}S_xF_{w_2}=C_{w_1w_0w_2^{-1}}S_x$ for $w_1,w_2\in B_0,x\in P^+$.\qed
\end{thm}
\textit{Proof.} It is easy to see $C_{z\pi}=C_z\ty{\pi}$ for any $z\in W'$, $\pi\in \Omega$. Using this fact, we can multiply some elements in $\Omega$ to $w_1w_0p_xw_2^{-1}$ to make the theorem \ref{uXi} applicable. Then we get $C_{w_1w_0p_xw_2^{-1}}=E_{w_1}C_{w_0p_x}F_{w_2}$. By \ref{low_c_0}(e), $C_{w_0p_x}=C_{w_0}S_x$, and hence the first equality follows. Since $S_x$  is in the center of $\mathcal{H}$, the second equality follows.\qed

The following  result is due to Guilhot  (see \cite{Guilhot2008lowest}).

\begin{cor}\label{left}
The left cells in the lowest two-sided cell of $\mathbf{c}_0$ of an extended affine Weyl group with positive weight function have the form $\Theta_w=\{w'w_vw^{-1}\vert w'\in U_v\}$, $w\in B_v$. In particular, the number of the left cells in $\mathbf{c}_0$ is equal to $|B_v|=|W_0|$.

\end{cor}
\textit{Proof.} Since we have known that  $\Theta_1=\{w'w_v\vert w'\in U_v\}$ is a left cell (see Lemma \ref{predecomp}(v)),  the corollary can be deduced directly from Theorem \ref{J0}.\qed

\begin{cor}
The  two-sided ideal $\mathcal{H}^{\mathbf{c}_0}=\mathbb{Z}[\Gamma]\{C_w \mid w\in\mathbf{c}_0\}$ is independent of the choice of the total order $\leq$ on $\Gamma$ and the positive weight function $L$.
\end{cor}
\textit{Proof.} By Theorem \ref{0xi}, $\mathcal{H}^{\mathbf{c}_0}$ is the (two-sided) ideal generated by $C_{w_0}$. Then the corollary follows from the formula $C_{w_0}=q_{w_0}^{-1}\sum_{w\in W_0}T_w$, see  Lemma \ref{predecomp}(i).
\begin{remark}

The  two-sided ideal  $\mathcal{H}^{\mathbf{c}_0}$ has affine cellularity \cite[Theorem 4.9]{Guilhot2014Cellularity}, in the sense of  \cite{koenig2012affine}. Indeed, the commutative ring in the definition of affine cell ideal is taken to be the center $\mathcal{Z}$. And the anti-involution is taken to be ${\cdot}^\flat$ (see the proof of Lemma \ref{predecomp}(iv)), noting that $(C_{w_1w_0p_xw_2^{-1}})^\flat=C_{w_2p_{-x}w_0w_1^{-1}}=C_{w_2w_0p_{x^{*}}w_1^{-1}}$ where $x^*=-w_0(x)$ and hence $(S_{x})^\flat=S_{x^*}$. Therefore, by Theorem \ref{0xi} one can conclude that $(C_{w_1w_0p_xw_2^{-1}},\flat)$ is a cellular basis of the lowest two-sided ideal. The cellular basis in \cite{Guilhot2014Cellularity} is actually taken to be $(C_{w_1w_0w_2^{-1}}P_x,\flat)$ where $P_x=S_{\omega_1}^{x_1}S_{\omega_2}^{x_2}\cdots S_{\omega_r}^{x_r}$ with $x=x_1\omega_1+\cdots+x_r\omega_r$ and $\omega_i$'s be the fundamental dominant weights. Therefore, the decomposition number studied in \cite[Sect. 6]{Guilhot2014Cellularity} is related to the decomposition of tensor product of fundamental representations.
\end{remark}
\section{Lusztig's properties P1-P15}\label{conj}
Recall that $W$ is an extended affine Weyl group, $\Gamma$ is a totally ordered abelian group, and $L$ is a positive weight function of $W$. An element $\alpha\in\mathbb{Z}[\Gamma]$ can be written in the form $\sum_{i}a_iq^{\gamma_i}$ with $a_i\in \mathbb{Z}$, $\gamma_i\in\Gamma$. The degree map $\deg:\mathbb{Z}[\Gamma]\longrightarrow\Gamma $ is defined by $\deg(\sum_{i}a_iq^{\gamma_i})=\max\{\gamma_i~|~a_i\neq0 \}$. Lusztig's $\mathbf{a}$-function is define by $\textbf{a}(z):=\sup\{\deg(h_{x,y,z})~|~x,y\in W \}$. Write $h_{x,y,z}=\gamma_{x,y,z^{-1}}q^{\textbf{a}(z)}+$ terms of degrees less than $\textbf{a}(z)$, if $\textbf{a}(z)\neq \infty$. Let $e$ be the neutral element of $W$. Define $\Delta(z):=-\deg(\pp{e,z})$. Then $\pp{e,z}$ has the form $$\pp{e,z}=q^{-L(z)}+\cdots+n_zq^{-\Delta(z)}$$ where $L(z)\geq\Delta(z)$ and $n_z$ is a nonzero integer.
The $\mathbb{Z}[\Gamma]$-linear map $\tau:\mathcal{H}\longrightarrow\mathbb{Z}[\Gamma]$ is defined by $\tau(\ty{w})=\delta_{w,e}$. It is well known that
\begin{equation}\label{tau}
\tau(\ty{x}\ty{y})=\delta_{xy,e}.
\end{equation} The subset $\d$ of $\mathbf{c}_0$ is defined to be $$\d=\{z\in \mathbf{c}_0\mid\textbf{a}(z)=\Delta(z) \}.$$ The following proposition says that Lusztig's properties P1-P15 (formulated in \cite[14.2]{lusztig2003hecke}) hold for the lowest generalzied two-sided cell $\mathbf{c}_0$.
\begin{prop}\label{lus conj}Keep the assumptions above.
\item[P$1'$] For any $z\in \mathbf{c}_0$, we have $\Delta(z)\geq \textbf{a}(z)$.
\item[P$2'$] If $d\in \d$, $x,y\in \mathbf{c}_0$ satisfying $\gamma_{x,y,d}\neq 0$, then $x=y^{-1}$.
\item[P$3'$] For any $y\in \mathbf{c}_0$, there exists a unique $d\in \d$ such that  $\gamma_{y^{-1},y,d}\neq0$.
\item[P$4'$] $\mathbf{c}_0=\{z\in W~|~\textbf{a}(z)=L(w_0) \}$.
\item[P$5'$] If $d\in\d$, $y\in \mathbf{c}_0$ and $\gamma_{y^{-1},y,d}\neq 0$ then $\gamma_{y^{-1},y,d}=n_d=1$.
\item[P$6'$] If $d\in\d$, then $d^2=e$.
\item[P$7'$] If $x,y,z\in \mathbf{c}_0$, then $\gamma_{x,y,z}=\gamma_{y,z,x}$.
\item[P$8'$] If $x,y,z\in \mathbf{c}_0$ are such that $\gamma_{x,y,z}\neq0$, then $x\sim_L y^{-1}$, $y\sim_L z^{-1}$, $z\sim_L x^{-1}$.
\item[P$13'$] For any left cell $\Theta$ in $\mathbf{c}_0$, $\Theta$ contains a unique $d\in\d$, and  $\gamma_{x^{-1},x,d}\neq0$ for any $x\in \Theta$.
\item[P$15'$]\footnote{The notation here is taken from \cite{geck2011rep}.} In $\mathbb{Z}[\Gamma]\otimes_\mathbb{Z}\mathbb{Z}[\Gamma]$, we have $\sum_{y'\in \mathbf{c}_0}h_{x,y',y}\otimes h_{w,x',y'}=\sum_{y'\in \mathbf{c}_0}h_{x,w,y'}\otimes h_{y',x',y}$ for $w\in \mathbf{c}_0$ and $x,x'\in W$.
\end{prop}
\textit{Proof.} P$4'$ is just the description in \ref{low_c_0}(a), see \cite[6.13]{bremke1997generalized}. Now, for any $w\in B_0,w'\in U_0$ consider
\begin{align*}
\tau(C_{ww_0w'^{-1}})&=\tau(E_wC_{w_0}F_{w'})\\
&=\sum_{u,u',y}\pp{uw_0,ww_0}\pp{y,w_0}\pp{u'w_0,w'w_0}\tau(\ty{u}\ty{y}\ty{u'^{-1}})
\end{align*}
where $u$ runs over $u\leq ww_0$, $l(uw_0)=l(u)+l(w_0)$, $u'$  runs over $u'\leq w'w_0$, $l(u'w_0)=l(u')+l(w_0)$, and $y$ runs over $y\leq w_0$.
If $\tau(\ty{u}\ty{y}\ty{u'^{-1}})\neq0$, then $\tau(\ty{u}\ty{yu'^{-1}})\neq0$ which implies  that $u^{-1}=yu'^{-1}$ (see \eqref{tau}) and hence $y=e$, $u=u'$. Therefore
\begin{align*}
\pp{e,ww_0w'^{-1}}&=\tau(C_{ww_0w'^{-1}})\\
&=\sum_u\pp{uw_0,ww_0}q^{-L(w_0)}\pp{uw_0,w'w_0}\\
&=\begin{cases}
q^{-L(w_0)} \text{+ lower  degree terms} \qquad  \qquad\qquad\qquad\qquad\text{if }w=w'\\
\text{an element of $\mathbb{Z}[\Gamma]$ with degrees } <L(w_0) \quad\text{ if } w\neq w'.
\end{cases}
\end{align*}

Therefore $\Delta(z)\geq L(w_0)=\textbf{a}(z)$ for any $z\in \mathbf{c}_0$ and the set $\d$ of $z$ such that $\textbf{a}(z)=\Delta(z)$ is $\{ww_0w^{-1}\mid w\in B_0 \}$. Thus  P$1'$ and  P$6'$ are proved. The elements in $\d$  will be called distinguished involutions. We also see that $n_d=1$ for each $d\in\d$.

By Corollary \ref{left},
\begin{itemize}
\item[(a)] $\Theta_w:=\{w'w_0w~|~w'\in U_0\},w\in B_0$ are all the left cells in $\mathbf{c}_0$, and hence every left cell contains a unique element in $\d$.
\end{itemize}
For $y\in \mathbf{c}_0$, denote by $\Theta_y$ the left cell containing $y$. The unique distinguished involution in $\Theta_y$ is denoted by $d_y$. Now for $x,y\in \mathbf{c}_0$ consider the equation $$\tau(C_xC_y)=\sum_{z\in\Theta_y}h_{x,y,z}\pp{1,z}.$$
By (P1') and (a), for any $z\in \Theta_y$ such that $z\neq d_y$, we have $\Delta(z)>\textbf{a}(z)$, and hence, for such $z$, $h_{x,y,z}\pp{1,z}$ has negative degree. Thus, $$\tau(C_xC_y)=\sum_{z\in\Theta_y}h_{x,y,z}\pp{1,z}\equiv h_{x,y,d_y}\pp{1,d_y}\equiv \gamma_{x,y,d_{y}}\mod{\mathbb{Z}[\Gamma^{<0}]}.$$
On the other hand, by the definition of KL-basis, we see $\tau(C_xC_y)\equiv \delta_{xy,1}\mod{\mathbb{Z}[\Gamma^{<0}]}$. Hence, $\gamma_{x,y,d_{y}}= \delta_{xy,1}$. Obviously, $\gamma_{x,y,d}=0$ if $d\neq d_y$ since $h_{x,y,d}=0$.

Parts P$2'$, P$3'$, P$5'$ and P$13'$ are proved.

P$7'$ follows from the two facts below.
\begin{itemize}
\item[(b)] For $z\in \mathbf{c}_0$, $\gamma_{x,y,z}$ equals the coefficient of $q_{w_0}$ in $\tau(\ty{x}\ty{y}\ty{z})$, see \cite[Cor.4.5]{bremke1997generalized} and the references therein.
\item[(c)] $\tau(hh')=\tau(h'h)$ for any $h,h'\in \mathcal{H}$.
\end{itemize}
P$8'$ follows from P$7'$, see\cite[14.8]{lusztig2003hecke}.

P$15'$ follows from Lemma \ref{comm} below.
\qed

Let $\mathcal{E}$ be a free $\mathbb{Z}[\Gamma]\otimes \mathbb{Z}[\Gamma]$ module with basis $\{\mathcal{E}_w\vert w\in \mathbf{c}_0\}$. We define a left module structure of $\mathcal{H}$ on $\mathcal{E}$ by
\[
C_x\mathcal{E}_w=\sum_{z\in W}( h_{x,w,z}\otimes 1)\mathcal{E}_z\quad \text{ for } x\in W \text{ and } w\in \mathbf{c}_0
\]
and define a right module structure of $\mathcal{H}$ on $\mathcal{E}$ by
\[
\mathcal{E}_w C_y=\sum_{z\in W} (1\otimes h_{w,y,z})\mathcal{E}_z\quad \text{ for } y\in W \text{ and } w\in \mathbf{c}_0.
\]
We will write $h_{x,y,z}=h_{x,y,z}\otimes1$ and $h'_{x,y,z}=1\otimes h_{x,y,z}$.
\begin{lemma}\label{comm} $\mathcal{E}$ is an $\mathcal{H}$-bimodule with above actions, i.e.
the left and right module structures are commutative. In particular, P$15'$ holds.
\end{lemma}

\textit{Proof. }The following claim is needed.

(g) Let $u$, $u'$ be both in $U_v$, then $(C_{uw_v}\mathcal{E}_{w_v})C_{w_vu'^{-1}}=C_{uw_v}(\mathcal{E}_{w_v}C_{w_vu'^{-1}})$.

Since $C_{uw_v}C_{w_v}=E_uC_{w_v}C_{w_v}=h_{w_v,w_v,w_v}E_uC_{w_v}=h_{w_v,w_v,w_v}C_{uw_v}$, we have $C_{uw_v}\mathcal{E}_{w_v}=h_{w_v,w_v,w_v}\mathcal{E}_{uw_v}$.

Similarly $C_{uw_v}C_{w_vu'^{-1}}=h_{w_v,w_v,w_v}E_uC_{w_v}F_{u'}$. Thus, $E_uC_{w_v}F_{u'}$ is $\overline{\,\cdot\,}$-invariant, hence we can write $E_uC_{w_v}F_{u'}=\sum_{z\in \mathbf{c}_0}b_zC_z$ with $b_{z}$ in $\mathbb{Z}[\Gamma]$ and $\overline{\,\cdot\,}$-invariant. Thus $h_{uw_v,w_vu'^{-1},z}=b_zh_{w_v,w_v,w_v}$, $z\in \mathbf{c}_0$. Since $\mathbf{a}(z)=L(w_0)=\deg(h_{w_v,w_v,w_v})$, $b_z=\gamma_{uw_v,w_vu'^{-1},z^{-1}}\in \mathbb{Z}$. Hence $\mathcal{E}_{uw_v}C_{w_vu'^{-1}}=h'_{w_v,w_v,w_v} \sum_{z\in \mathbf{c}_0}\gamma_{uw_v,w_vu'^{-1},z}\mathcal{E}_z$. Then

\begin{align*}
(C_{uw_v}\mathcal{E}_{w_v})C_{w_u'^{-1}}&=h_{w_v,w_v,w_v}\mathcal{E}_{uw_v}C_{w_vu'^{-1}}\\
&=h_{w_v,w_v,w_v}h'_{w_v,w_v,w_v} \sum_{z\in \mathbf{c}_0}\gamma_{uw_v,w_vu'^{-1},z}\mathcal{E}_z
\end{align*}
Similar computations show that $C_{uw_v}(\mathcal{E}_{w_v}C_{w_vu'^{-1}})=h_{w_v,w_v,w_v}h'_{w_v,w_v,w_v} \sum_{z\in \mathbf{c}_0}\gamma_{uw_v,w_vu'^{-1},z}\mathcal{E}_z$, and thus the claim (g) is proved.

Now we prove the lemma. For any $x,y\in W$, $w\in B_v$, $w'\in U_v$, we have
\begin{align*}
(C_x\mathcal{E}_{ww_vw'^{-1}})C_y&=\frac1{h_{w_v,w_v,w_v}}(C_x(C_{ww_v}\mathcal{E}_{w_vw'^{-1}}))C_y\\
&=\frac1{h_{w_v,w_v,w_v}}((C_xC_{ww_v})\mathcal{E}_{w_vw'^{-1}})C_y\\ \intertext{(We can write $C_xC_{ww_v}=\sum_{u\in U_v}h_{x,ww_v,uw_v}C_{uw_v}$ since $\{uw_v|u\in U_v \}$ is a left cell.)}
&=\frac1{h_{w_v,w_v,w_v}}(\sum_{u\in U_v}h_{x,ww_v,uw_v}C_{uw_v}\mathcal{E}_{w_vw'^{-1}})C_y \\
&=\frac1{h_{w_v,w_v,w_v}h'_{w_v,w_v,w_v}} (\sum_{u\in U_v}h_{x,ww_v,uw_v}C_{uw_v}(\mathcal{E}_{w_v}C_{w_vw'^{-1}}))C_y
\intertext{(by claim (g) we have)}\\
&=\frac1{h_{w_v,w_v,w_v}h'_{w_v,w_v,w_v}} (\sum_{u\in U_v}h_{x,ww_v,uw_v}(C_{uw_v}\mathcal{E}_{w_v})C_{w_vw'^{-1}})C_y\\
&=\frac{h_{w_v,w_v,w_v}}{h_{w_v,w_v,w_v}h'_{w_v,w_v,w_v}} (\sum_{u\in U_v}h_{x,ww_v,uw_v}(\mathcal{E}_{uw_v}C_{w_vw'^{-1}}))C_y\\
&=\frac{h_{w_v,w_v,w_v}}{h_{w_v,w_v,w_v}h'_{w_v,w_v,w_v}} \sum_{u\in U_v}h_{x,ww_v,uw_v}\mathcal{E}_{uw_v}(\sum_{u'\in U_v}h'_{w_vw'^{-1},y,w_vu'^{-1}}C_{w_vu'^{-1}})\\
&=\frac1{h_{w_v,w_v,w_v}h'_{w_v,w_v,w_v}} \sum_{\substack{u\in U_v\\u'\in U_v}}h_{x,ww_v,uw_v}h'_{w_vw'^{-1},y,w_vu'^{-1}}(C_{uw_v}\mathcal{E}_{w_v})C_{w_vu'^{-1}}
\end{align*}

Similarly we can get \[C_x(\mathcal{E}_{ww_vw'^{-1}}C_y)=\frac1{h_{w_v,w_v,w_v}h'_{w_v,w_v,w_v}} \sum_{\substack{u\in U_v\\u'\in U_v}}h_{x,ww_v,uw_v}h'_{w_vw'^{-1},y,w_vu'^{-1}}C_{ww_v}(\mathcal{E}_{w_v}C_{w_vu'^{-1}}).\] Using claim (g) again, we see that $(C_x\mathcal{E}_{ww_vw'^{-1}})C_{y}=C_x(\mathcal{E}_{ww_vw'^{-1}}C_y)$. Now the lemma follows since each $z\in \mathbf{c}_0$ has the form $z=ww_vw'^{-1}$ for some $w\in B_v$, $w'\in U_v$. \qed

Now taking the coefficients of $1\otimes q_{w_{0}}$ in the equality in P$15'$, one obtains the following Corollary:
\begin{cor}\label{cor_comm}
For any $x,z\in W$, $y\in \mathbf{c}_0$, we have \[
\sum_{w\in \mathbf{c}_0,v\in \mathbf{c}_0}h_{x,y,w}\gamma_{w.z.v^{-1}}=\sum_{v\in \mathbf{c}_0,w\in \mathbf{c}_0}h_{w,z,v^{-1}}\gamma_{y,z,w^{-1}}.
\]\qed
\end{cor}
\section{Based rings of $\mathbf{c}_0$}\label{base}

Let  $J_0$ be the free $\mathbb{Z}$-module  with basis $\{t_x\mid x\in \mathbf{c}_0 \}$. Define multiplication in $J_0$ by\begin{equation}\label{mul}
t_xt_y=\sum_{z\in\textbf{c}_0}\gamma_{x,y,z^{-1}}t_z.
\end{equation}
where the integer $\gamma_{x,y,z^{-1}}$ is the coefficient of $q^{\textbf{a}(z)}$ in $h_{x,y,z}$. According to \cite[\S18]{lusztig2003hecke} we see that Prop. \ref{lus conj} implies the following result.

\begin{thm_def}\label{defJ}
With the multiplication \eqref{mul}, the $\mathbb{Z}$-module $J_0$ becomes an associative ring with  with identity $\sum_{d\in \d}t_d$. The ring $J_0$ is called the {\em based ring} of the lowest generalized two-sided cell $\textbf{c}_0$.
\end{thm_def}

Moreover, according to \cite[Thm. 18.9]{lusztig2003hecke}, we have

\begin{thm}\label{maptoJ}
The $\mathbb{Z}[\Gamma]$-linear map $
\phi :\mathcal{H}\longrightarrow \mathbb{Z}[\Gamma]\otimes_\mathbb{Z}J_0$ defined by
\[ \quad C_x\mapsto \sum_{d\in\d,z\in \mathbf{c}_0}h_{x,d,z}t_z\]is  a homomorphism of $\mathbb{Z}[\Gamma]$-algebras  with 1.
\end{thm}

\begin{thm}\label{based ring} If $u=w_1w_0p_xw_2^{-1}$, $u'=w_3w_0p_{x'}w_4^{-1}$, $u''=w_5w_0p_{x''}w_6^{-1}$, where $w_i\in B_0$, $i=1,\dots,6$ and $x,x',x''\in P^+$, then $\gamma_{u,u',u''^{-1}}\neq0$ implies that $w_2=w_3$,  $w_4=w_6$, $w_5=w_1$, and $\gamma_{u,u',u''^{-1}}=m(x,x',x'')$ where $m(x,x',x'')$ is the multiplicity of the simple module $V(x'')$ in the tensor product $V(x)\otimes V(x')$.
\end{thm}
\textit{Proof.} We follow the approach of \cite[Cor.4.3]{Xi1990based}.

By P$8'$ in Proposition \ref{lus conj}, $\gamma_{u,u',u''^{-1}}\neq0$ implies that $u \sim_L u'^{-1}$, $u'\sim_L u''$, and $u''^{-1}\sim_L u^{-1}$. Then by Corollary \ref{left} about the description of the left cells in $\textbf{c}_0$, we have $w_2=w_3$, $w_4=w_6$ and $w_5=w_1$. By Lemma \ref{predecomp} and property \ref{low_c_0} (e), 
\begin{align*}
C_{u}C_{u'}&=C_{w_1w_0p_xw_2^{-1}}C_{w_2w_0p_{x'}w_4^{-1}}\notag\\&=E_{w_1p_{w_0x}}C_{w_0w_2^{-1}}C_{w_2w_0}F_{(p_{x'}w_4^{-1})^{-1}}\\
&=\sum_{y\in P^+}h_{w_0w_2^{-1},w_2w_0,w_0p_y}E_{w_1p_{w_0x}}C_{w_0p_y}F_{(p_{x'}w_4^{-1})^{-1}}\\
&=\frac1{h_{w_0,w_0,w_0}^2}\sum_{y\in P^+}h_{w_0w_2^{-1},w_2w_0,w_0p_y}C_{w_1w_0p_{x}}C_{w_0p_y}C_{w_0p_{x'}w_4^{-1}}\\
&=\sum_{x'',y\in P^+}h_{w_0w_2^{-1},w_2w_0,w_0p_y}a_{y,x''}C_{w_1w_0p_{x''}w_4^{-1}}
\end{align*}
for some integer $a_{y.x''}$. Furthermore, we have \begin{equation}\label{g1}
a_{0,x''}=m(x,x',x'').
\end{equation}
Thus  we have  $h_{u,u',w_1w_0p_{x''}w_4^{-1}}=\sum_{y\in P^+}a_{y,x''}h_{w_0w_2^{-1},w_2w_0,w_0p_y}$. Hence, we have
\begin{equation}\label{g2}
\gamma_{u,u',(w_1w_0p_{x''}w_4^{-1})^{-1}}=\sum_{y\in P^+}a_{y,x''}\gamma_{w_0w_2^{-1},w_2w_0,(w_0p_y)^{-1}}
\end{equation}
Then the theorem would follows from \eqref{g1}, \eqref{g2}, and the following claim:
\begin{equation}\label{conf}
\text{If }y\neq0,\text{ then }\gamma_{w_0w_2^{-1},w_2w_0,(w_0p_y)^{-1}}=0.
\end{equation}
By P$7'$, we have $\gamma_{w_0w_2^{-1},w_2w_0,(w_0p_y)^{-1}}=\gamma_{(w_0p_y)^{-1},w_0w_2^{-1},w_2w_0}$. Thus to prove \eqref{conf}, we only need to prove $h_{(w_0p_y)^{-1},w_0w_2^{-1},w_0w_2^{-1}}$ has degree $<L(w_0)$ if $y\neq0$. But this follows immediately from the equality $C_{p_{-y}w_0}C_{w_0w_2^{-1}}=C_{p_{-y}w_0w_2^{-1}}$ (see Theorem \ref{0xi}). This completes the proof. \qed
\begin{cor}\label{iso}
Let $\operatorname{Mat}_{B_0}(\z)$ be the matrix ring  over $\z$  and indexed by $B_0$, see Definition \ref{center form} for the definition of $\z$. Then there is an isomorphism of rings between the based ring $J_0$ and the matrix ring $\operatorname{Mat}_{B_0}(\z)$ given by
 $$J_0\longrightarrow \operatorname{Mat}_{B_0}(\z), t_{w_1w_0p_xw_2^{-1}}\mapsto S_xI_{w_1,w_2}$$ where $I_{w_1,w_2}$ is the matrix whose $(w_1,w_2)$-entry is 1 and  other entries are 0.
\end{cor}
\textit{Proof.} This is follows from the last theorem and the formula $$S_{x}S_{x'}=\sum_{x''\in P^+}m(x,x',x'')S_{x''}$$ (see Section \ref{low_c_0} (d)).\qed
\begin{remark}

From Theorem \ref{based ring}, we see that $\gamma_{x,y,z}$, $x,y,z\in \mathbf{c}_0$ is actually independent of the choice of the positive weight function,  and so is the based ring $J_0$. There is one way to see it directly. Since $\gamma_{x,y,z}$ is the coefficients of $q_{w_0}$ in $m_{x,y,z^{-1}}:=\tau(\ty{x}\ty{y}\ty{z})$(see Section \ref{conj}(e)), it is independent of the choice of the total order on $\Gamma$.

Assume that $q_1=q_s,q_2=q_t$ with $s$ and $t\in S$ are not conjugate in $W$ (if it is possible). Write $\xi_1=q_1-q_1^{-1}$, $\xi_2=q_2-q_2^{-1}$ and denote $\nu_1$ (resp. $\nu_2$) the number of simple reflections conjugate to $s$ (resp. $t$) in a reduced expression of $w_0$. Obviously, $\nu_1+\nu_2=l(w_0)$.  Then we can write $m_{x,y,z}$ in the form \begin{gather*}
m_{x,y,z}=\sum_{\substack{0\leq i\leq \nu_1\\0\leq j\leq\nu_2}}a_{i,j}\xi_1^i\xi_2^j
\end{gather*}with $a_{i.j}\in\mathbb{Z}$; the proof is similar to that of  \cite[Thm. 7.2]{Lusztig1985cells}. Then we can conclude that $\gamma_{x,y,z}=a_{\nu_1,\nu_2}$ is independent of the choice of the positive weight function.

In this way, we also see the based ring $J_0$ is well defined for any positive weight function, since it is well know that $J_0$ is well defined when the weight function is constant. But  to construct the homomorphism $\phi$ in Theorem \ref{maptoJ}, it seems that the Proposition \ref{lus conj} is unavoidable.\qed
\end{remark}

\section{Type $\tilde C_n$ case}\label{nonext}

In this section the affine Weyl group $W'$ is of type $\tilde C_r$ $(r\ge 1)$.   Here we identify $\tilde{A}_1$ with $\tilde{C}_1$, and $\tilde{B}_2$ with $\tilde{C}_2$. We will discuss the based ring of lowest generalized two-sided cell in this case for completeness. Note that in this case, some weight functions on $W'$ can not be extended to $W$.

\begin{figure}
\centering
\includegraphics[width=0.6\linewidth]{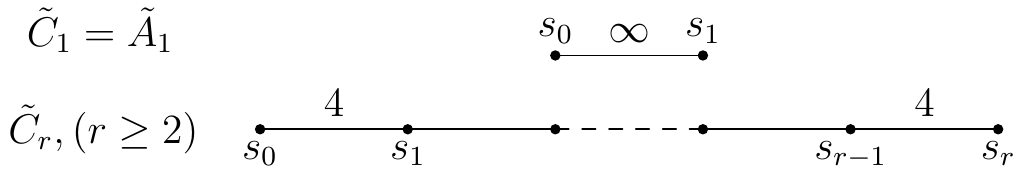}
\caption{Dynkin diagram of type $\tilde{C}_r,r\geq1$.}
\label{fig:diag-eps-converted-to}
\end{figure}

We number the simple reflections of $W'$ as usual, see Fig. \ref{fig:diag-eps-converted-to}.
Let $L$ be a positive weight function on $W'$ such that $L(s_0)\neq L(s_r)$. It is no harm to assume that $L(s_0)< L(s_r)$.


Now the claim in \ref{extended aff}(a) does not hold in this case.
We have to modify the definition of special points. A point $\lambda\in E$ is called special if $m(\lambda):=\sum_{\lambda\in H}L(H)$ is maximal, see \cite{bremke1997generalized,bremke1996Phd_thesis}. It turns out that the set of special points is $Q=\mathbb{Z}R$. Define the quarters $\mathcal{C}_\lambda^+,\lambda\in Q$, as in subsection 2.3.  The boxes are defined to be the connected components of the complement  in $E$ of the set of hyperplanes supporting the walls of the quarters $\mathcal{C}_\lambda^+,\lambda\in Q$.
The box containing $\lambda(\in Q)$ in its closure and contained in $\mathcal{C}_\lambda^+$ is still denoted by $\Pi_\lambda^+$. Now each box contains $|W_0|$ alcoves (while in the case $L(s_0)=L(s_r)$ each box contains $|W_0|/|\Omega|=|W_0|/2$ alcoves.)

Then the lowest generalized two sided cell of $(W',L)$ is \[
\mathbf{c}'_0=\{w_1w_0p_xw_2^{-1}\,|\,w_1,w_2\in B_0,x\in Q^+ \}
\]
where $B_0=\{w\in W'\,|\,wA_0^+\subset\Pi_0^+ \}$, $w_0$ is the longest element $W_0$, and $Q^+=Q\cap P^+$. It can also be characterized by the $\textbf{a}$-value\[
\mathbf{c}'_0=\{z\in W'\,|\,\textbf{a}(z)=L(w_0) \}.
\]
One can see \cite{bremke1997generalized,bremke1996Phd_thesis} for these results.

The discussions in the Section \ref{formula}--\ref{base} are also applicable in current situation. We state it here without proof. Since Lemma \ref{lus_lemma} still holds (see \cite[Lemma 4.3]{bremke1997generalized}), we have decomposition (compare with Theorem \ref{J0},\ref{uXi})\[
C_{w_1w_0p_xw_2^{-1}}=E_{w_1}C_{w_0p_x}F_{w_2}=E_{w_1}C_{w_0}F_{(p_xw_2^{-1})^{-1}}
\]
for any $w_1,w_2\in B_0$, $x\in Q^+$. Then the Lusztig conjectures on cells (P1-P15) hold for $\mathbf{c}'_0$. The proof is completely similar to Proposition \ref{lus conj}. Then we have the based ring $J'_0=\mathbb{Z}\{t_w\,|\,w\in \mathbf{c}'_0 \}$ with structure constant $\gamma_{x,y,z^{-1}}$ and the homomorphism from $\mathcal{H'}$ to $\mathbb{Z}[\Gamma]\otimes_\mathbb{Z}J'_0$, see Theorem \ref{maptoJ}.

For any $x\in Q$, write $x=y-z$ with some $y,z\in Q^+$ and define $\theta_x$ to be $\ty{y}(\ty{z})^{-1}$. Also define, for $x\in Q^+$,
$$S_x=\sum_{x'\in Q}d(x',x)\theta_{x'}$$
where $d(x',x)$ is the dimension of the weight space $V(x)_{x'}$ as Section \ref{low_c_0}. Then $S_x,x\in Q^+$ is a $\mathbb{Z}[\Gamma]$-basis of the center of $\mathcal{H}'$. Actually this holds for any affine Hecke algebra $\mathcal{H'}$.  It is easy to see that if $x\in Q^+$, then $d(x',x)\neq0$ only if $x'\in Q$. Then one can see, for $x,x'\in Q^+$,
\[
S_{x}S_{x'}=\sum_{x''\in Q^+}m(x,x',x'')S_{x''},
\]
where $m(x,x',x'')$ is multiplicity of $V(x'')$ in the tensor product $V(x)\otimes V(x')$. Following  the same lines as \cite{lusztig1983singularities} or \cite[2.9]{xi1994representations}, we have
$$C_{w_0p_x}=C_{w_0}S_{x},\text{ for } x\in Q^+.$$
And hence  $C_{w_0p_x}C_{w_0p_y}=\sum_{z\in Q^+}m(x,y,z)h_{w_0,w_0,w_0}C_{w_0z}$ for $x,y\in Q^+$. Then the structure of the based ring $J_0$ can be explicitly described as
\[
t_{w_1w_0p_xw_2^{-1}}t_{w_3w_0p_yw_4^{-1}}=\sum_{z\in Q^+}\delta_{w_2,w_3}m({x,y,z})t_{w_1w_0p_zw_4^{-1}}
\]
where $w_1,w_2,w_3,w_4\in B_0$, $x,y\in Q^+$. In other words, we have

\medskip

{\bf Proposition 6.2.} The based ring $J'_0$ of the lowest generalized two-sided cell $\mathbf{c}'_0$ of $W'$  is isomorphic to a matrix ring of order $|W_0|$ over the commutative ring $\mathcal{Z'}_{\mathbb{Z}}$, where $\mathcal{Z'}_{\mathbb{Z}}=\mathbb{Z}\langle S_x\,|\,x\in Q^+ \rangle$ (which is isomorphic to a subring of  the representation ring of the algebraic group $G_\mathbb{C}$).

\medskip

{\bf Remark 6.3.} If $L$ is a positive weight function of $W'$ such that $L(s_0)=L(s_r)$, then $L$ can be extended to a weight function on $W$. In this case the lowest generalized two-sided cell $\mathbf{c}''_0$ of $W'$ equals $W'\cap \mathbf{c}_0$, where $\mathbf{c}_0$ is the lowest generalized two-sided cell of $W$. Note that $\mathbf{c}''_0$ is different from $\mathbf{c}'_0$ (see subsection 6.1 for the definition of $\mathbf{c}'_0$) and the based ring of $\mathbf{c}''_0$ is a subring of $J_0$ but it  is not isomorphic to the based ring $J'_0$ of $\mathbf{c}'_0$.

\section{Irreducible representations attached to $\mathbf{c}_0$}\label{rep1}

 Let $k$ be a field. Assume that there is a ring homomorphism $\mathbb{Z}[\Gamma]\longrightarrow k$.  Set $\mathcal{H}_k:=\mathcal{H}\otimes_{\mathbb{Z}[\Gamma]}k$, $J_{0,k}:=J_{0}\otimes_\mathbb{Z}k$. The homomorphism $\mathcal{H}_k\longrightarrow J_{0,k}$ induced by $\phi$ in Theorem \ref{based ring} is  denoted again by $\phi$. Denote by $\mathcal{H}_k^{\mathbf{c}_0}$ the  two-sided ideal of $\mathcal{H}_k$ spanned by all $C_w\otimes 1$, $w\in\textbf{c}_0$. We shall simply write  $C_w$ and $t_w$ instead of  $C_w\otimes 1$  and $t_w\otimes 1$ respectively.
 
 In this section we use the explicit structure of $J_0$ to study representations of $\mathcal{H}_k$ attached to $\mathbf{c}_0$. We shall follow the approaches in \cite{lusztig1987cellsIII,xi2007representations}. Given a $J_{0,k}$-module $E$, through the homomorphism $\phi:\mathcal{H}_k\longrightarrow J_{0,k}$, we get an  $\mathcal{H}_k$-module structure on $E$. This $\mathcal{H}_k$-module will be denoted by $E_\phi$ or $\phi_*(E)$. We will see that oftern $E_{\phi}$ has a unique simple quotient when $E$ is irreducible, see Lemma \ref{unique_max_submod} below.
 
 \begin{defn}
 We say that a  representation $M$ of $\mathcal{H}_k$  is attached to $\mathbf{c}_0$ if $H_k^{\mathbf{c}_0}M\neq0$. The set of all irreducible representations of $\mathcal{H}_k$  attached to $\mathbf{c}_0$ is denoted by $\irr(\mathcal{H}_k,\mathbf{c}_0)$. The set of  all irreducible representations of $\mathcal{H}_k$  is denoted by $\irr(\mathcal{H}_k)$.

\end{defn}
Since $\mathcal{H}_k^{\mathbf{c}_0}$ is a two-sided ideal generated by $C_{w_0}$ (see Theorem \ref{0xi}), we have
\begin{itemize}
\item[(a)] $\irr(\mathcal{H}_k, \mathbf{c}_0)=\{M\in \irr \mathcal{H}_k~|~C_{w_0}M\neq0 \}$.
\end{itemize}

Note that $J_{0,k}$ is naturally a left $J_{0,k}$-module by multiplication in $J_{0,k}$.
Define a right $\mathcal{H}_{k}$-module structure on $J_{0,k}$ by \[t_wC_x=\sum_{v\in \mathbf{c}_0 }h_{w,x,v}t_v\quad \text{ for }w\in \mathbf{c}_0, x\in W.\] Then by Corollary \ref{cor_comm}
\begin{itemize}
\item[(b)]$J_{0,K}$ is a $J_{0,k}$-$\mathcal{H}_{k}$-bimodule.
\end{itemize}  Via the homomorphism $\phi :\mathcal{H}_{k}\longrightarrow J_{0,k}$, $J_{0,K}$ becomes an $\mathcal{H}_k$-bimodule. We have
\begin{itemize}
\item[(c)]$\mathcal{H}_k$ acts on $J_{0,K}$ by $C_xt_w=\sum_{v\in \mathbf{c}_0}h_{x,w,v}Z_v$ for $x\in W, w\in \mathbf{c}_0$,
\end{itemize}
In fact, $C_xt_w=\phi(C_x)t_w=\sum_{d\in\d,u\in \mathbf{c}_0}h_{x,d,u}t_ut_w=\sum_{\substack{d\in\d,\\u,v\in \mathbf{c}_0}}h_{x,d,u}\gamma_{u,w,v^{-1}}t_v$,
 which is $\sum_{v\in \mathbf{c}_0}h_{x,w,v}Z_v$ by $ P15' $, \(P2'\), \(P3'\), \(P5'\), \(P7'\).

Let $M$ be an $\mathcal{H}_k$-module. Then $\hat {M}:={J_{0,K}}\otimes_{\mathcal{H}_k} M$ is a  $J_{0,K}$-module.   Define $\tilde M=\phi_*(\hat M)$. It is easy to verify that
\begin{itemize}
\item[(d)]The map $\tilde M\longrightarrow M$, $t_w\otimes\ m\mapsto C_u m$ is a homomorphism of $\mathcal{H}_k$-modules. If $M$ is a simple $\mathcal H_k$-module with $C_{w_0}M\neq0$, then $M$ is the unique composition factor of $\tilde M$ such that $C_{w_0}M\neq0$
\end{itemize}

Let $E$ be a $J_{0,k}$-module and $N$ be an $\mathcal{H}_k$-submodule of $E_\phi$. Using Corollary \ref{cor_comm} we can see the following
\begin{itemize}
\item[(e)]$\widehat{N}\longrightarrow E$, $Z_u\otimes n\mapsto \phi(C_u).n$ is a homomorphism of $J_{0,k}$-modules, where $n\in N$ is viewed as an element of $E$.
\end{itemize}The image of this map is $\mathcal{H}_k^{\mathbf{c}_0}N$.

The following lemma is crucial for main results in this section, see \cite[Lemma 2.4]{xi2007representations} for the original version.
\begin{lemma}\label{unique_max_submod}
Assume that $E$ is a simple $J_{0,k}$-module such that $C_{w_0}E_\phi\neq 0$. Then the subspace $K=\{b\in E_\phi |C_u.b=0 \text{ for all }u\in \mathbf{c}_0\}$ is the unique maximal submodule of $E_\phi$.
In particular, $E_\phi$ has only one composition factor $M'$ such that $C_{w_0}M'\neq0$. And $M'$ is the unique simple quotient of $E_\phi$.\end{lemma}
\textit{Proof}. It is easy to see that $K$ is an $\mathcal{H}_k$-submodule of $E_\phi$.

Let $v\in E_\phi$ be an element not in $K$ and $N=\mathcal{H}_k. v$.  Then the image $\mathcal{H}_k^{\mathbf{c}_0}N$ of the homomorphism $\hat{N}\longrightarrow E$ is nonzero. Since $E$ is a simple $J_{0,k}$-module, we have $E=\mathcal{H}_k^{\mathbf{c}_0}N$. Thus $N=E_\phi$. Therefore $K$ is the unique maximal submodule of $E_\phi$. Other statements follow immediately.\qed

Denote by $\irr(J_{0,k},C_{w_0})$ the set of all irreducible representations $E$ of $ J_{0,k}$ such that $C_{w_0}E_{\phi}\neq0$. The above lemma says that there is a well defined map $$\rho:\irr(J_{0,k},C_{w_0})\longrightarrow\irr(\mathcal{H}_k,\mathbf{c}_0)$$ such that $\rho(E)$ is the unique simple quotient of $E_\phi$.
\begin{prop}\label{irr_bijection}
$\rho$ is bijective.
\end{prop}
\textit{Proof.} First we prove that $\rho$ is surjective. Let $M\in \irr(\mathcal{H}_k,\mathbf{c}_0)$. Then the map $\tilde M\longrightarrow M$ has nonzero image $\mathcal{H}_k^{\mathbf{c}_0}M\neq0$ since $C_{w_0}M\neq0$, and hence surjective. Since $\tilde M=\phi_*(\hat M)$, $\hat M$ must have a composition factor $E$ such that $E_\phi$ has a composition factor $M$. Since $C_{w_0}M\neq0$, by Lemma \ref{unique_max_submod}, $M$ is the unique simple quotient of $E$, i.e. $\rho(E)=M$. Therefore $\rho$ is surjective.

Now we prove that $\rho$ is injective. Let $E\in\irr(J_{0,k},C_{w_0})$ and $\pi:E_\phi\longrightarrow M=\rho(E)$ be the quotient map. Let $p':\hat{E}\longrightarrow \hat M$, $t_u\otimes x\mapsto t_u\otimes \pi(x)$ be the  homomorphism of $J_{0,k}$-modules induced from $\pi$. Then we have commutative diagram
\[
\xymatrix{\hat{E}\ar[r]^{p'}\ar[d]_{\theta} & \hat M\ar[d]^p\\E_\phi\ar[r]^\pi& M}
\]
where $\theta$ is the homomorphism of $J_{0,k}$-module defined in (e) and  $p$ is essentially the map  defined in (d).

The map $p'$ induces a surjective homomorphism $\overline{p'}:\hat{E}/\ker\theta\longrightarrow \hat M/p'(\ker\theta)$ of $J_{0,k}$-modules. On one hand, $C_{w_0}M\neq0$ implies that $p$ is surjective and hence $\ker p\neq \hat M$. On the other hand, the commutative diagram implies that $p'(\ker \theta)\subset \ker p$. Then we have $\hat M/p'(\ker\theta)\neq0$. Since $\widehat{E}/\ker\theta\simeq E$ is simple, $\overline{p'}$ is an isomorphism. Therefore $E$ is a composition factor of $\hat M$.

We claim that $\hat M$ admits one and only one composition factor $E'$ such that $C_{w_0}E'_\phi\neq0$. Otherwise, by Lemma \ref{unique_max_submod}, $\tilde M$ has more than one composition factor $M'$ such that $C_{w_0}M'\neq0$, which contradicts  (d). Since $C_{w_0}E_\phi\neq0$, $E$ is unique characterized as the simple factor of $\hat M$ on which the action of $C_{w_0}$ is nonzero. Thus $\rho$ is injective. This completes the proof.\qed

\medskip

By  proposition above, to parametrize the irreducible representations of $\mathcal{H}_k$ attached to $\mathbf{c}_0$, it suffices to parametrize the irreducible representations $E$ of $J_{0,k}$ such that $C_{w_0}E_{\phi}\neq0$.


\def\mspec{\text{Max}(\mathcal{Z}_k)}

The set of simple modules (up to isomosphism) of $\mathcal{Z}_k=\mathcal Z_{\mathbb Z}\otimes k$  are one-to-one correspondence to the set $\mspec$ of maximal ideals of  $\mathcal{Z}_k$. For any $\m\in\mspec$, the field $k_\m:=\mathcal{Z}_k/\m$ is a simple $\mathcal{Z}_k$-module. Denote by $\lambda_\m$ by the quotient map $\mathcal{Z}_k\rightarrow k_\m$.

Recall that $J_0$ is isomorphic to $\operatorname{Mat}_{B_0}(\z)$ via $t_{w_1w_0p_xw_2^{-1}}\mapsto S_xI_{w_1,w_2}$ where $I_{w_1,w_2}$ is the matrix whose $(w_1,w_2)$-entry is 1 and other entries are 0.
Thus, given a maximal ideal $\m\in\mspec$, the homomorphism $J_{0,k}\longrightarrow \operatorname{Mat}_{B_0}(k_\m)$, $t_{w_1w_0p_xw_2^{-1}}\mapsto \lambda_\m(S_x)I_{w_1,w_2}$ gives a representations of $J_{0,k}$. The affording space will be denoted by $E_\m$.
Then
\begin{itemize}
\item[(f)]The map $\m\mapsto E_\m$ gives a bijection between the set $\mspec$ and the set of simple representations of $J_{0,k}$ over $k$. And $\dim_{k} E_\m=(\dim_{k_\m} E_\m)[k_\m:k]=|W_0|[k_\m:k]$.
\end{itemize}
We denote $(E_\m)_\phi$ simply  by $E_{\m,\phi}$.
\begin{prop}\label{condition1}
$C_{w_0}E_{\m,\phi}=0$ if and only if $\sum_{x\in P^+}h_{w_0,ww_0,w_0p_x}S_x\in\m$ for all $w\in B_0$. In particular, $C_{w_0}E_{\m,\phi}\neq0$ for any maximal ideal $\m\in\mspec$ if $h_{w_0,w_0,w_0}$ is nonzero in $k$.
\end{prop}
\textit{Proof.} Note that $$\phi(C_{w_0})=\sum_{w\in B_0,x\in P^+}h_{w_0,ww_0w^{-1},w_0p_xw^{-1}}t_{w_0p_xw^{-1}}=\sum_{w\in B_0,x\in P^+}h_{w_0,ww_0,w_0p_x}t_{w_0p_xw^{-1}}.$$

Thus $C_{w_0}$ acts on $E_{\m,\phi}$ by the matrix $$\sum_{w\in B_0,x\in P^{+}}h_{w_0,ww_0,w_0p_x}\lambda_\m(S_x)I_{e,w}.$$
Hence, $C_{w_0}E_{\m,\phi}=0$ if and only if, for any $w\in B_0$, $\sum_{x\in P^+}h_{w_0,ww_0,w_0p_x}\lambda_\m(S_x)=0$, i.e. $\sum_{x\in P^+}h_{w_0,ww_0,w_0p_x}S_x\in\m$. The second statement follows from the fact that when $w=e$, $\sum_{x\in P^+}h_{w_0,ww_0,w_0p_x}S_x=h_{w_0,w_0,w_0}T_e$ where $T_e$ is the identity of $\mathcal{H}$.  \qed

Here are some notations. $I_0=S\cap W_0$. For any subset $I\subset I_0$, the dominant weight $x_I$ is defined to be $\sum_{i\in I}\omega_i$, where $\omega_i$ is the fundamental dominant weight corresponding to $i$. Let $w_I$ denote the longest element in the parabolic subgroup $W_I$ of $W'$. And denote the complement of $I$ in $I_0$ by $I'$. Then (see the proof in \cite[3.7]{xi1994based})
\begin{itemize}
\item[(g)] $x_Iw_{I'}=ww_0$ for some $w\in B_0$.
\end{itemize}
Define \begin{align}\label{ad}
\alpha_I:=&\sum_{x\in P^+}h_{w_0,x_Iw_{I'},w_0p_x}S_x\in\mathcal{Z}_k:=\mathcal{Z}\otimes_{\mathbb{Z}[\Gamma]}k,\\
\zeta_I:=& h_{w_I,w_I,w_I}\in k,\\
\Delta_k:=&\{I\subset I_0~|~\zeta_{I'}\neq0,\text{ but }\zeta_{I'\cup\{i \}}=0 \text{ for any } i \in I \}.\label{ad3}
\end{align}

Then the following proposition is proved in \cite[3.7]{xi1994based} in one parameter case. One can easily generalize it to the case of positive weight function in view of the preparation in Section \ref{formula} and Section \ref{conj}.
\begin{prop}\label{condition2}If $\m\in\mspec$, then
$\sum_{x\in P^+}h_{w_0,ww_0,w_0p_x}S_x\in\m$ for all $w\in B_0$ if and only if $\alpha_I\in\m$ for all $I\in \Delta_k$.
\end{prop}

Now combining Proposition \ref{irr_bijection}, Proposition \ref{condition1} with Proposition \ref{condition2}, we obtain
\begin{thm}\label{basic set}
Fix a field $k$ and a specialization $\mathbb{Z}[\Gamma]\longrightarrow k$. The simple modules of $\mathcal{H}_k$ attached to $\mathbf{c}_0$ have a bijection with the set $\{\m\in\mspec\,|\,\alpha_I\notin\m \text{ for some }I\in \Delta_k\}$,  where $\alpha_I$ and $\Delta_k$ are defined by (\ref{ad}) (\ref{ad3}).\qed
\end{thm}

\begin{remark}
\item[(i)]Note that if $h_{w_0,w_0,w_0}\neq0$, then $\emptyset\in \Delta_k$ and $\alpha_{\emptyset}=h_{w_0,w_0,w_0}\notin\m$, for any $\m\in \mspec$. So in this case the simple $\mathcal{H}_k$-modules attached to $\mathbf{c}_0$ are in bijection with the set $\mspec$.
\item[(ii)] If $k$ is algenraically closed, then $\mspec$ has a bijection with the set of representatives of semisimple conjugacy classes in the simply connecte simple algebraic group $G_\mathbb{k}$ over $k$ with root system  $R$.
The  correspondence is as follows.
If $s\in G_\mathbb{C}$ is a semisimple element, then we can define a homomorphism $\lambda_s:\mathcal{Z}_k\ra k$, $\lambda_s(S_x)=tr(s,V(x))$ for $x\in P^+$, where $V(x)$ is a simple $G_\mathbb{C}$-module with highest weight $x$. Then $\ker \lambda_s\in\mspec$.
\item[(iii)]When $W'$ is of type $\tilde C_n$ using the results in section 6 we can get similar results for representations of Hecke algebras $\mathcal H'_k$ of $W'$.
\end{remark}

We conclude this section with a dimension formula.

For any $\mathcal{H}_k$-module $M$, denote by $\operatorname{Ann}_M\mathbf{c}_0$ the largest submodule of $M$ annihilated by the lowest two-sided ideal $\mathcal{H}_k^{\mathbf{c}_0}$. For any $J_{0,k}$-module $E$ we denote $$\rho(E)=E_\phi/\operatorname{Ann}_{E_\phi}\mathbf{c}_0.$$

\begin{prop}\label{kdim}
Let $\m\in\mspec$. Then the dimension of $\rho(E_\m)$ over $k_\m$ is the rank of the matrix $(\lambda_\m(m_{w,w'}))_{w,w'\in B_0}$, where \[m_{w,w'}=\sum_{x\in P^+}h_{w_0w^{-1},w'w_0,w_0p_x}S_x\in \mathcal{Z}_k.\]
\end{prop}
\textit{Proof}. By Theorem \ref{J0}, we have \begin{equation}\label{l1}
\operatorname{Ann}_{E_\m}\mathbf{c}_0=\{v\in E_{\m,\phi}\,\vert\, C_{w_0w^{-1}}v=0 \text{ for all } w\in B_0\}.
\end{equation}

One can see that $\phi(C_{w_0w^{-1}})= \sum_{\substack{w'\in B_0\\x\in P^+}}h_{w_0w^{-1},w'w_0w'^{-1},w_0p_xw'^{-1}}t_{w_0p_xw'^{-1}}$ acts on $E_\m$ by the matrix
\begin{align*}
&\sum_{\substack{w'\in B_0\\x\in P^+}}h_{w_0w^{-1},w'w_0w'^{-1},w_0p_xw'^{-1}}\lambda_\m(S_x)I_{e,w'}\\
&=\sum_{\substack{w'\in B_0\\x\in P^+}}h_{w_0w^{-1},w'w_0,w_0p_x}\lambda_\m(S_x)I_{e,w'}\\
&=\sum_{w'\in B_0}\lambda_\m(m_{w,w'})I_{e,w'}.
\end{align*}Thus, by \eqref{l1},
we have $\dim_{k_\m}\operatorname{Ann}_{E_\m}\mathbf{c}_0=\# B_0-\operatorname{rank}(\lambda(m_{w,w'}))_{w,w'\in B_0}$. Therefore \[\dim_{k_\m}(\rho(E_\m))=\operatorname{rank}(\lambda_\m(m_{w,w'}))_{w,w'\in B_0}.\]\qed

\section{Representations over the residue field of the center}\label{rep2}

In this section, we provide another perspective for the representations of affine Hecke algebras. Recall that $\h$ is an affine Hecke algebra over $\mathbb{Z}[\Gamma]$, $\mathcal{Z}$ is the center of $\h$, $J_0$ is the based ring of the lowest generalized two-sided cell $\mathbf{c}_0$. Denote $\jg=J_0\otimes_\mathbb{Z}\mathbb{Z}[\Gamma]$. We will see that the map $\phi:\h\ra\jg$ is injective. And we will find out the set of prime ideals of $\mathcal{Z}$ such that $\phi_\pr:k_\pr\h\ra k_\pr\jg$ is an isomorphism, where $k_\pr$ is the residue field of the center $\mathcal{Z}$ at $\pr$. In particular, $\pr=0$ is such an ideal.
\begin{lemma}
Denote $\Delta=\h C_{w_0}$.
\begin{itemize}
\item[(i)] $\Delta$ is an $\h$-$\mathcal{Z}$-bimodule.
\item[(ii)] $\Delta$ is a faithful left $\h$-module.
\item[(iii)] $\Delta$ is a free $\mathcal{Z}$-module
 with basis $\{C_{ww_0}~|~w\in B_0 \}$.
 \end{itemize}
\end{lemma}
\textit{Proof.} (i) is obvious.

The proof of (ii) is similar to \cite[1.7]{lusztig1987cellsIII}. First as claimed in \cite[1.1(c)]{lusztig1987cellsIII}, we have
\begin{itemize}
\item[(a)] For any $y\in W$ we can find $s_1,\cdots, s_p\in S$ such that $l(ys_1s_2\cdots s_p)=l(y)+p$ and $l(ys_1s_2\cdots s_pt)>l(ys_1s_2\cdots s_p)$ for any $t\in S_0:=S\cap W_0$.
\end{itemize}
Now assume $0\neq h\in\h$. Write $h=\sum a_w\ty{w}$. Let $y$ be a maximal $w$ (under Bruhat order) such  that $a_w\neq0$. Then one can find $s_1,s_2,\cdots,s_p$ as in (a). Let $h'=h\ty{s_1}\cdots\ty{s_p}$. It is easy to see that the coefficient of $\ty{ys_1\cdots s_p}$ in $h'C_{w_0}$ is $a_y$, which is nonzero by hypothesis. Thus $h'C_{w_0}\neq0$. Hence we have proved that there exists $u\in \Delta=\h C_{w_0}$ such that $h.u\neq0$. Therefore, $\mathcal{H}C_{w_0}$ is a faithful $\h$-module.

(iii). First, by Lemma \ref{predecomp}(v) or by Theorem \ref{J0}, $\Delta$ is a $\mathbb{Z}[\Gamma]$-module generated by $\{C_{ww_0}|w\in U_0 \}$. By Theorem \ref{0xi}, for each $w\in U_0$, there exists $w'\in B_0$, $x\in P^+$ such that $C_{ww_0}=C_{w'w_0p_x}=C_{w'w_0}S_x$. Thus $\Delta$ is generated by $\{C_{w'w_0}|w'\in B_0 \}$ over $\mathcal{Z}$. Second, $\{C_{w'w_0}|w'\in B_0 \}$ is free over $\mathcal{Z}$ since  $\{C_{ww_0}|w\in U_0 \}$ is free over $\mathbb{Z}[\Gamma]$.
\qed

\medskip

The $\h$-$\mathcal{Z}$-bimodule structure on $\Delta$ naturally induces a $\zg$-algebra homomorphism \begin{equation}\label{map1}
\varphi:\h\ra\operatorname{End}_\mathcal{Z}(\Delta).
\end{equation}And by Corollary \ref{iso}, there is an isomorphism of $\zg$-algebras
\begin{align}
\eta : &\jg\ra \operatorname{End}_\mathcal{Z}(\Delta),\label{map2}\\& t_{w_1w_0p_xw_2^{-1}}\mapsto (C_{w_3w_0}\mapsto \delta_{w_2,w_3}S_xC_{w_1w_0}).
\end{align}
\text{where $w_1,w_2,w_3\in B_0$}.

\begin{thm}\label{commutative}

The following diagram is commutative,
\[
\xymatrix{\h\ar[r]^{\phi}\ar[rd]_{\varphi} & \jg\ar[d]^\eta_\simeq\\{}& \ezd}
\]
where the map $\phi$ is defined in Theorem \ref{maptoJ}, the map $\varphi$ is defined in \eqref{map1}, and the map $\eta$ is defined in \eqref{map2}.
\end{thm}
\textit{Proof.} Since $\{C_{ww_0}\,|\,w\in B_0\}$ is a $\mathcal{Z}$-basis of $\Delta$, it suffices to prove that for any $w\in W$, $w_1\in B_0$, $(\eta\circ\phi)(C_w)(C_{w_1w_0})=\varphi(C_w)(C_{w_1w_0})=C_w C_{w_1w_0}$.
\begin{align*}
LHS&=\eta(\sum_{w_3,w_2\in B_0,x\in P^{+}}h_{w,w_3w_0w_3^{-1},w_2w_0p_xw_3^{-1}}t_{w_2w_0p_xw_3^{-1}})(C_{w_1w_0})\\
&=\sum_{w_2,x}h_{w,w_1w_0w_1^{-1},w_2w_0p_xw_1^{-1}}C_{w_2w_0p_x}\\
&=\sum_{w_2,x}h_{w,w_1w_0,w_2w_0p_x}C_{w_2w_0p_x}=RHS
\end{align*}\qed

From this theorem, one can see that $\phi$ and $\varphi$ are the same if one identifies $J_0$ with $\ezd$.

\begin{cor}\label{y1}
\begin{itemize}
\item[(i)] $\phi:\h\ra \jg$ is injective.
\item[(ii)]$\phi(\mathcal{Z})$ is the center of $\jg$.
\end{itemize}
\end{cor}
\textit{Proof.} Since $\Delta$ is a faithful $\h$-module, $\varphi$ is injective. Thus $\phi$ is injective and (i) follows. Since the center of $\operatorname{End}_\mathcal{Z}(\Delta)$ is $\varphi(\mathcal{Z})$, the center of $\jg$ is $\phi(\mathcal{Z})$. So (ii) follows.\qed

\medskip

For any prime ideal $\mathfrak{p}\in\operatorname{Spec}\mathcal{Z}$, the residue field of $\mathcal{Z}$ at $\pr$ is denoted by $k_\pr$, which by definition is the field of fractions of the quotient $\mathcal{Z}/\pr$. The natural morphism $\mathcal{Z}\ra k_\pr$ is denoted by $\lambda_\pr$. Set
\begin{align*}
k_\pr\mathcal{H}:=k_\pr\otimes_\mathcal{Z}\mathcal{H}, \quad k_\pr\jg&:=k_\pr\otimes_\mathcal{Z}\jg\quad k_\pr\Delta=k_\pr\otimes_\mathcal{Z}\Delta\\
\phi_\pr&:k_\pr\mathcal{H}\ra k_\pr\jg\\
m_{w_1,w_2}&:=\sum_{x\in P^+}h_{w_0w_1^{-1},w_2w_0,w_0p_x}S_x\in\mathcal{Z}\\
\det&:=\det(m_{w_1,w_2})_{w_1,w_2\in B_0}\in\mathcal{Z}.
\end{align*}
Both $\mathcal{H}$ and $\jg$ are free $\mathcal{Z}$-module with rank $|W_0|^2$, see \cite[2.10]{xi1994representations}. Thus $\dim_{k_\pr}k_\pr\mathcal{H}=\dim_{k_\pr}k_\pr\jg=|W_0|^2$.
\begin{lemma}\label{y4}
$\phi_\pr$ is an isomorphism if and only if $k_\pr\Delta$ is a simple module of $k_\pr\mathcal{H}$.
\end{lemma}
\textit{Proof.} By Theorem \ref{commutative}, $k_\pr\mathcal{H}$ acts on $k_\pr\Delta$ via the homomorphism $\phi_{\pr}:k_\pr\mathcal{H}\ra k_\pr\jg\simeq \operatorname{End}_{k_\pr}(k_\pr \Delta)$. If $k_\pr\Delta$ is a simple $k_\pr\mathcal{H}$-module (which must be of dimension $|W_0|$), then by Artin-Wedderburn theorem, $\phi_\pr(k_\pr\mathcal{H})$ must be of dimension $\geq|W_0|^2$. Thus $\phi_\pr$ is surjective,  and hence bijective. The another direction is easy.\qed

\begin{prop}\label{y2}
$\phi_\pr$ is an isomorphism if and only if $\det\notin\pr$.
\end{prop}
\textit{Proof.} First assume that  $\phi_\pr$ is an isomorphism.  Then $k_\pr\mathcal{H}$ is a simple algebra over $ k_\pr $,  and $k_\pr\Delta$ is its only  simple module up to isomorphism. Let $ h $ be the image in the algebra $k_\pr\mathcal{H}$. Since $k_\pr\Delta=(k_\pr\mathcal{H})h  $ is a simple module, we have $ h\neq 0 $. Moreover, $ h^2\neq 0 $; otherwise, $ h $ is a nonzero nilpotent element in $k_\pr\mathcal{H}$, which contradicts with that $k_\pr\mathcal{H}$ is a simple algebra. Hence $h(k_\pr\Delta)\neq0$. Then the quotient of $k_\pr\Delta$ by the submodule consisting of the elements annihilated by the ideal $ (k_\pr\mathcal{H})h(k_\pr\mathcal{H}) $ is nonzero (see Lemma \ref{unique_max_submod}) and has the dimension $\operatorname{rank}(\lambda_\pr(m_{w_1,w_2}))_{w_1,w_2\in B_0}$ (see Proposition \ref{kdim}). But $k_\pr\Delta$ is a simple module. We see that the matrix $ (\lambda_\pr(m_{w_1,w_2}))_{w_1,w_2\in B_0} $ is of full rank and hence $\lambda_\pr(\det)\neq0$, i.e. $\det\notin\pr$.

Conversely, if $\lambda_\pr(\det)\neq0$, then  $k_\pr\Delta$ is a simple module (see Proposition \ref{kdim}), and hence $\phi_\pr$ is an isomorphism by Lemma \ref{y4}.\qed

\begin{thm}\label{last}
Let $K$ be the field of fractions of $\mathcal{Z}$, and let $\mathcal{Z}[\frac1\det]$ be the subring of $K$ generated by $\mathcal{Z}$ and $\frac1\det$. Then $\phi:\mathcal{H}[\frac1\det]\ra\jg[\frac1\det]$ is an isomorphism of $\mathcal{Z}[\frac{1}{\det}]$-algebras. In particular, $K\mathcal{H}\simeq  K\jg$ is a split simple algebra over $K$.
\end{thm}
\textit{Proof.}  By Corollary \ref{y1}, $\phi:\mathcal{H}[\frac1\det]\ra\jg[\frac1\det]$ is injective. 

Note that the prime ideals $ \mathfrak{q} $ of $ \mathcal{Z}[\frac1\det] $  is in bijection with the prime ideals $ \pr $ of $ \mathcal{Z} $ with $ \det\notin \pr $. Precisely, $ \mathfrak{q}=\mathfrak{p}[\frac1\det] $.
Since $ (\mathcal{Z}[\frac1\det])_\mathfrak{q}=\mathcal{Z}_\pr $, we know  that the localization $ (\mathcal{H}[\frac1\det])_\mathfrak{q}\ra(\jg[\frac1\det])_\mathfrak{q} $  is an isomorphism for any prime ideal  $ \mathfrak{q} $ of $ \mathcal{Z}[\frac1\det] $ (see Proposition \ref{y2}). Then we conclude that  $ \phi:\mathcal{H}[\frac1\det]\ra\jg[\frac1\det] $ is an isomorphism by the following well-known fact about commutative algebras (applying to $ A=\mathcal{Z}[\frac1\det] $, $ M=\mathcal{H}[\frac1\det] $, and $ N= \jg[\frac1\det]$):
\begin{itemize}
\item[]Let $A$ be a noetherian commutative ring, and $f:M\ra N$ be a morphism of finitely generated $A$-modules. If $f_\mathfrak{q}:k_\mathfrak{q} M\ra k_\mathfrak{q} N$ is surjective for any prime ideal $\mathfrak{q}$ of $A$, then $f$ is surjective.
\end{itemize}\qed

\textbf{Acknowledgment:} This paper is a part of PhD thesis of the author under the supervision of Professor Nanhua Xi. The author thanks the guide and helpful comments of Professor Nanhua Xi. This work is partly supported by the National Natural Science Foundation of China (No.11321101).

\bibliography{reff}

\end{document}